\documentclass{article}
\usepackage{amsmath,amsthm,amssymb}

\newtheorem{thm}{Theorem}[section]
\newtheorem{crl}[thm]{Corollary}

\newtheorem{lmm}[thm]{Lemma}

\newtheorem{prp}[thm]{Proposition}

\theoremstyle{definition}
\newtheorem{dfn}[thm]{Definition}

\theoremstyle{remark}
\newtheorem*{rem}{Remark}

\def\H{\mathfrak{H}}
\def\K{\mathfrak{K}}
\def\Z{\mathbb{Z}}
\def\N{\mathbb{N}}
\def\RHS{\mathrm{RHS}}
\def\Out{\mathrm{Out}}
\def\RQC{\mathrm{RQC}}
\def\Fix{\mathrm{Fix}}
\def\Stab{\mathrm{Stab}}
\def\SL{\mathrm{SL}}
\def\ne{\mathrm{ne}}

\title{The universal relatively hyperbolic structure on a group
and relative quasiconvexity for subgroups}
\date{}
\author{
Yoshifumi Matsuda\footnote
{Graduate School of Mathematical Sciences, 
University of Tokyo, 
3-8-1 Komaba, 
Meguro-ku, 
Tokyo, 
153-8914 Japan,
ymatsuda@ms.u-tokyo.ac.jp} \footnote{
The author is supported by the Global COE Program at Graduate School of Mathematical Sciences, the University of Tokyo, and Grant-in-Aid for Scientific Researches for Young Scientists (B) (No. 22740034), Japan Society of Promotion of Science.}, 
Shin-ichi Oguni\footnote
{Department of Mathematics, Faculty of Science,
Ehime University,
2-5 Bunkyo-cho, 
Matsuyama, 
Ehime, 
790-8577 Japan,
oguni@math.sci.ehime-u.ac.jp}, 
Saeko Yamagata\footnote
{Faculty of Education and Human Sciences, Yokohama National University,
240-8501 Yokohama, Japan,
yamagata@ynu.ac.jp}
}

\begin{document}

\maketitle
\begin{abstract}     
We discuss the notion of the universal relatively hyperbolic structure on a group 
which is used in order to characterize relatively hyperbolic structures on the group.
We also study relations between relatively hyperbolic structures on a group 
and relative quasiconvexity for subgroups of the group.\\


\noindent
2010MSC:
20F67;
20F65
\end{abstract}

\section{Introduction}\label{intro}
The notion of relatively hyperbolic groups was introduced in \cite{Gro87} 
and has been studied by many authors 
(see for example \cite{Bow12}, \cite{D-S05b}, \cite{Far98} and \cite{Osi06a}). 
When a countable group is relatively hyperbolic, 
relative quasiconvexity for subgroups can be defined. 
There are several equivalent definitions of relative hyperbolicity for countable groups 
and those of relative quasiconvexity for subgroups (see \cite[Section 3 and Section 6]{Hru10}). 
In this paper we adopt definitions in terms of geometrically finite convergence actions. 

We regard a conjugacy invariant collection of subgroups of a countable group 
relative to which the group is hyperbolic as a structure on the group 
and call such a collection a relatively hyperbolic structure on the group. 
In this paper we discuss the notion of 
the universal relatively hyperbolic structure on a group 
which is used in order to characterize relatively hyperbolic structures on the group.
In particular we give a characterization of 
the universal relatively hyperbolic structure on a finitely generated group
(Corollary \ref{char-universal}). 
We also study relations between relatively hyperbolic structures on a group 
and relative quasiconvexity for subgroups of the group. 
Indeed we give two theorems: 
\begin{itemize}
\item (Theorem \ref{inf-rhs}) when a countable group is not virtually cyclic 
and admits a proper relatively hyperbolic structure, 
it has two families of infinitely many relatively hyperbolic structures 
such that structures in one family have pairwise distinct collections of 
relatively quasiconvex subgroups while
structures in the other family have the same collection of 
relatively quasiconvex subgroups; 
\item (Theorem \ref{qc-updown2}) two relatively hyperbolic structures on 
a finitely generated group have the same collection of relatively quasiconvex subgroups 
if and only if these structures are equal when we ignore virtually infinite cyclic subgroups.
\end{itemize}

The following are contents of this paper. 
Section 2 gives preliminaries. We recall several facts about convergence actions, 
and give definitions of relatively hyperbolic structures on a countable group 
and relatively quasiconvex subgroups. 
In Section 3, blow-ups and blow-downs of relatively hyperbolic structures are defined. 
By using these notions, 
we define a partial order on the set of all relatively hyperbolic structures on a countable group. 
In Section 4, the universal relatively hyperbolic structure is defined in order to characterize relatively hyperbolic structure on a countable group. 
In Section 5, we consider how relative quasiconvexity for subgroups varies when we blow up and down a relatively hyperbolic structure.
In Section 6, cardinality of the set of relatively hyperbolic structures is studied and 
Theorem \ref{inf-rhs} is proved. 
In section 7, when we consider a finitely generated group, 
it is shown that the partially ordered set of all relatively hyperbolic structures 
on the group is a directed set (Proposition \ref{fiberprod}). 
Also Corollary \ref{char-universal} is proved. 
In section 8, we show Theorem \ref{qc-updown2}. 
In Appendix A, examples of torsion-free countable groups without the universal relatively hyperbolic structure are given. 
In Appendix B, we discuss a condition for a characterization of relatively hyperbolic structures by using mapping class  groups. 

Here we introduce several notations. 
In this paper 
$G$ denote a countable group with the discrete topology. 
Let $L$ be a subgroup of $G$. 
Let $\H$ and $\K$ be two conjugacy invariant collections of infinite subgroups of $G$.
We define $L$-conjugacy invariant collections 
$\H_L$ and $L \wedge \H$ of infinite subgroups of $L$, 
and $G$-conjugacy invariant collections 
$\H_{\ne}$ and $\H\wedge \K$ of infinite subgroups of $G$ as 
\begin{align*}
\H_L=&\{ H\in \H ~|~ H \subset L\};\\
L \wedge \H=&\left\{ P ~\left|~ 
P\text{ is an infinite subgroup of }L, 
P=L \cap H\text{ for some }H \in \H
\right.
\right\};\\
\H_{\ne}=&\H \setminus \{ H \in \H ~|~ H\text{ is virtually infinite cyclic}\};\\
\H\wedge \K=&\left\{ P ~\left|~ 
P\text{ is an infinite subgroup of }G, 
P=H \cap K\text{ for some }H\in \H\text{ and }K\in \K  \right.
\right\}.
\end{align*}

\section{Relatively hyperbolic structures on a group 
and relative quasiconvexity for subgroups}\label{Relatively}

First we recall some definitions and properties related to convergence actions (refer to \cite{G-M87}, \cite{Tuk94}, \cite{Tuk98}, \cite{Bow12} and \cite{Bow99b}). 
Let $G$ have a continuous action on a compact metrizable space $X$. 
The action is called a convergence action if $X$ has infinitely many points and for each infinite sequence $\{g_i\}$ of mutually different elements of $G$, there exist a subsequence $\{g_{i_j}\}$ of $\{g_i\}$ and two points $r,a\in X$ such that $g_{i_j}|_{X\setminus \{r\}}$ converges to $a$ uniformly on each compact subset of $X\setminus \{r\}$ 
and also $g_{i_j}^{-1}|_{X\setminus \{a\}}$ converges to $r$ uniformly on each compact subset of $X\setminus \{a\}$. 
The sequence $\{g_{i_j}\}$ is called a convergence sequence and also the points $r$ and $a$ are called the repelling point of $\{g_{i_j}\}$ and the attracting point of $\{g_{i_j}\}$, respectively. 

We fix a convergence action of $G$ on a compact metrizable space $X$. 
The set of all repelling points and attracting points is equal to the limit set $\Lambda(G,X)$ (\cite[Lemma 2M]{Tuk94}). 
The cardinality of $\Lambda(G,X)$ is $0$, $1$, $2$ or $\infty$ (\cite[Theorem 2S, Theorem 2T]{Tuk94}). 
If $\#\Lambda(G,X)=\infty$, then the action of $G$ on $X$ is called a 
non-elementary convergence action and also the induced action of $G$ on $\Lambda(G,X)$ is a minimal non-elementary convergence action.
We remark that $\#\Lambda(G,X)=0$ if $G$ is finite by definition. 
Also $G$ is virtually infinite cyclic if $\#\Lambda(G,X)=2$ (see \cite[Lemma 2Q,Lemma 2N and Theorem 2I]{Tuk94}). 
An element $l$ of $G$ is said to be loxodromic 
if it is of infinite order and has exactly two fixed points. 
For a loxodromic element $l\in G$, the sequence $\{l^i\}_{i\in \N}$ is a convergence sequence with the repelling point $r$ and the attracting point $a$, 
which are distinct and fixed by $l$. 
We call a subgroup $H$ of $G$ a parabolic subgroup if 
it is infinite, fixes exactly one point and has no loxodromic elements. 
Such a point is called a parabolic point. 
A parabolic point is said to be bounded if its maximal parabolic subgroup 
acts cocompactly on its complement. 
We call a point $r$ of $X$ a conical limit point if there exists a convergence sequence $\{g_i\}$ with the attracting point $a\in X$ such that the sequence $\{g_i(r)\}$ converges to a different point from $a$.  
The convergence action is said to be geometrically finite 
if every point of $X$ is either a conical limit point or a bounded parabolic point.
Since $X$ has infinitely many points, every geometrically finite convergence action is non-elementary. 
Also since all conical limit points and all bounded parabolic points belong to the limit set, every geometrically finite convergence action is minimal.

Based on \cite[Definition 1]{Bow12} and \cite[Theorem 0.1]{Yam04} 
(see also \cite[Definition 3.1]{Hru10}), 
we define relatively hyperbolic structures on a countable group from a dynamical viewpoint.
\begin{dfn}\label{def-rhs}
Let $\H$ be a conjugacy invariant collection of infinite subgroups of $G$. 
The group $G$ is said to be hyperbolic relative to $\H$ 
if there exists a convergence action of $G$ on a compact metrizable space 
satisfying the following: 
\begin{enumerate}
\setlength{\itemsep}{0mm}
\item[(1)] the set of all maximal parabolic subgroups of the action is equal to $\H$; 
\item[(2)] one of the following holds:
    \begin{enumerate}
    \setlength{\itemsep}{0mm}
    \item[(i)]  the limit set of the action is a finite set; 
    \item[(ii)] the limit set of the action is an infinite set and the induced convergence action of $G$ on the limit set is geometrically finite. 
    \end{enumerate}
\end{enumerate}
Such a collection $\H$ is called a relatively hyperbolic structure on $G$. 
We remark that any relatively hyperbolic structure on $G$ has only 
finitely many conjugacy classes (see \cite[Theorem 1B]{Tuk98}).
A relatively hyperbolic structure is said to be trivial (resp. proper) if it is (resp. is not) equal to $\{G\}$. 
We denote the set of all relatively hyperbolic structures on $G$ by $\RHS(G)$. 
\end{dfn}

Let $G$ be endowed with a relatively hyperbolic structure $\H$ and 
consider a convergence action of $G$ satisfying the conditions in Definition \ref{def-rhs}. 
It essentially follows from \cite[Theorem 9.4]{Bow12} 
and \cite[Theorem 0.1]{Yam04} that the limit set of such a convergence action 
is uniquely determined as a $G$-space independently of the choice of a convergence action. 
We denote such a $G$-space by $\partial (G,\H)$, which is called the Bowditch boundary. 
Indeed we have the following:
\begin{prp}\label{boundary}
Let $\H$ be a relatively hyperbolic structure on $G$. 
Consider a convergence action of $G$ on a compact metrizable space $X$ satisfying the following: 
\begin{enumerate}
\setlength{\itemsep}{0mm}
\item[(1)] the set of all maximal parabolic subgroups of the action is equal to $\H$;
\item[(2)] one of the following holds:
    \begin{enumerate}
    \setlength{\itemsep}{0mm}
    \item[(i)]  the limit set of the action is a finite set; 
    \item[(ii)] the limit set of the action is an infinite set and the induced convergence action of $G$ on the limit set is geometrically finite. 
    \end{enumerate}
\end{enumerate}
Then the limit set $\Lambda(G,X)$ of the action is uniquely determined 
as a $G$-space independently of the choice of a convergence action and a compact metrizable space $X$. 
\end{prp}
\begin{proof}
If $\Lambda(G,X)$ is an infinite set, then the assertion follows from \cite[Theorem 9.4]{Bow12} and \cite[Theorem 0.1]{Yam04}. 
Now we suppose that the limit set $\Lambda(G,X)$ is a finite set. 
Under this assumption, the proof is divided into the following three cases. 

In the case where $\Lambda(G,X)$ is empty, the group $G$ is finite and $\H$ is empty. 
Hence the limit set of every convergence action of $G$ satisfying the conditions above is empty. 

In the case where $\Lambda(G,X)$ consists of one point, $\H$ is equal to $\{G\}$. 
The limit set of every convergence action of $G$ satisfying the conditions above consists of one point. 

In the case where $\Lambda(G,X)$ consists of two points, it follows from \cite[Theorem 2Q and Theorem 2I]{Tuk94} that there exists a loxodromic element $l \in G$ such that the subgroup $\langle l \rangle$ generated by $l$ is of finite index. 
Replacing $l$ by its power if necessary, we may assume that $\langle l \rangle$ is a finite index normal subgroup of $G$. 
For every convergence action of $G$ 
such that the limit set consists of two points, 
$l$ is loxodromic and the limit set consists of the two fixed point of $l$. 
For each element $g \in G$, either $glg^{-1}=l$ or $glg^{-1}=l^{-1}$ holds. 
If the former holds, then $g$ acts trivially on the limit set. 
If the latter holds, then $g$ interchanges the two points. 
Thus the action of $G$ on the limit set is independent of a choice of the convergence action and a compact metrizable space $X$. 
\end{proof}

The following are special examples of relatively hyperbolic structures. 
\begin{itemize}
\setlength{\itemsep}{0mm}
\item $G$ is infinite if and only if it has the trivial relatively hyperbolic structure $\{G\}$. 
Moreover, $\partial(G, \{G\})$ consists of a single point. 
\item $G$ is hyperbolic if and only if 
the empty collection $\emptyset$ is a relatively hyperbolic structure on $G$. 
Moreover, $\partial(G, \emptyset)$ is the Gromov boundary of $G$ 
(see \cite[8.2]{Gro87} and also \cite[Theorem 3.4 and Theorem 3.7]{Fre95}). 
\end{itemize}

We can refine Proposition \ref{boundary} for the case of 
virtually infinite cyclic groups. 
We note that a virtually infinite cyclic group $G$ has the maximal finite normal subgroup $K$ and $G/K$ is isomorphic to either $\Z$ or $\Z/2\Z*\Z/2\Z$ (see for example \cite[Lemma 11.4]{Hem76}).
\begin{crl} 
Let $G$ be a virtually infinite cyclic group and $K$ be the maximal finite normal subgroup of $G$. 
Then $G/K$ is isomorphic to $\Z$ if and only if the action of $G$ on $\partial(G,\emptyset)$ is trivial. 
\end{crl}
\begin{proof}
We put $G_0=G/K$ and define a generating set $S_0$ of $G_0$ as follows. 
\begin{align*}
S_0=\left\{
\begin{array}{ll}
\{t\} & ~{\text{if}}~ G_0 ~{\text{is isomorphic to}}~ \Z \\
\{a, b\} & ~{\text{if}}~ G_0 ~{\text{is isomorphic to}}~ \Z/2\Z*\Z/2\Z,
\end{array}
\right. 
\end{align*}
where $a,b \in G_0$ satisfy $a^2=b^2=1$. 
The Cayley graph $\Gamma(G_0, S_0)$ is hyperbolic and the Gromov boundary $\partial\Gamma(G_0, S_0)$ consists of two points. 
We consider a compact metrizable space $\overline{\Gamma(G_0, S_0)}=\Gamma(G_0, S_0) \cup \partial\Gamma(G_0, S_0)$ and the action of $G$ on $\overline{\Gamma(G_0, S_0)}$ which is induced from the action of $G_0$ on $\Gamma(G_0, S_0)$. 
This action is a convergence action and the limit set $\Lambda(G, \overline{\Gamma(G_0, S_0)})$ is equal to $\partial\Gamma(G_0, S_0)$. 
If $G_0$ is isomorphic to $\Z$, then the action of $G$ on $\partial\Gamma(G_0, S_0)$ is trivial. 
If $G_0$ is isomorphic to $\Z/2\Z*\Z/2\Z$, then the element $a \in G_0$ acts on $\Gamma(G_0, S_0)$ as an inversion and hence the action of $G$ on $\partial\Gamma(G_0, S_0)$ is nontrivial. 
Since $\partial\Gamma(G_0, S_0)$ is $\partial(G,\emptyset)$ by Proposition \ref{boundary}, the assertion follows. 
\end{proof}

We recall the notion of relative quasiconvexity for subgroups of a countable group with a relatively hyperbolic structure in accordance with \cite[Definition 1.6]{Dah03}. 

\begin{dfn}\label{def-relqc}
Let $G$ be endowed with a relatively hyperbolic structure $\H$. 
A subgroup $L$ is said to be quasiconvex relative to $\H$ in $G$ if one of the following holds:
\begin{enumerate}
\setlength{\itemsep}{0mm}
\item[(1)] $\partial(G, \H)$ is a finite set; 
\item[(2)] $\partial(G, \H)$ is an infinite set and $\Lambda (L, \partial(G,\H))$ is a finite set; 
\item[(3)] Both $\partial(G, \H)$ and $\Lambda (L, \partial(G,\H))$ are infinite sets and the induced convergence action of $L$ on $\Lambda (L, \partial(G,\H))$ is geometrically finite. 
\end{enumerate}
We denote the set of all subgroups of $G$ 
that are quasiconvex relative to $\H$ in $G$ by $\RQC(G, \H)$. 
\end{dfn}

\begin{rem}\label{ele-relqc}
Let $G$ be endowed with a relatively hyperbolic structure $\H$.
\begin{enumerate}
\setlength{\itemsep}{0mm}
\item[(I)]  
It immediately follows from the definition that for every $H \in \H$, every subgroup of $H$ is quasiconvex relative to $\H$ in $G$. 
Also since the limit set of a convergence action of a virtually abelian group is a finite set (see \cite[Theorem 2U]{Tuk94}), every virtually abelian subgroup of $G$ is quasiconvex relative to $\H$ in $G$. 
In particular every finite subgroup of $G$ is quasiconvex relative to $\H$ in $G$. 
\item[(II)] When $\H=\emptyset$, 
a subgroup $L$ of $G$ is quasiconvex relative to $\emptyset$ in $G$ 
if and only if it is quasiconvex in the ordinary sense (see \cite[Proposition 4.3]{Bow99b}).
\item[(III)] If a subgroup $L$ of $G$ is quasiconvex relative to $\H$, then $L \wedge \H$ is a relatively hyperbolic structure on $L$ by \cite[Theorem 9.1]{Hru10}. 
\end{enumerate}
\end{rem}

\section{A partial order on the set of relatively hyperbolic structures}
We define a relation $\to$ on 
the set of conjugacy invariant collections of infinite subgroups of $G$ as follows. 
For two conjugacy invariant collections $\H$ and $\K$ of infinite subgroups of $G$, 
$\K \to \H$ holds if for every $K \in \K$, there exists $H \in \H$ such that $K \subset H$. 
When $\K \to \H$ holds, $\K$ is called a blow-up of $\H$ and $\H$ is called a blow-down of $\K$
(compare with \cite{M-O-Y5}). 
We show the following.  

\begin{prp}\label{order}
The relation $\to$ defines a partial order on $\RHS(G)$. 
\end{prp}
Note that the outer automorphism group $\Out(G)$ of $G$ acts naturally on $\RHS(G)$ and that this action preserves the order $\to$.

\begin{dfn}
Let $G$ be a group and let $\H$ be a conjugacy invariant collection of infinite subgroups of $G$. 
The collection $\H$ is said to be almost malnormal in $G$ if the following hold: 
\begin{enumerate}
\setlength{\itemsep}{0mm}
\item[(1)] the intersection of every pair of two elements of $\H$ is finite. 
\item[(2)] every element of $\H$ is equal to its normalizer in $G$. 
\end{enumerate}
\end{dfn}

A subgroup $H$ of $G$ is said to be almost malnormal in $G$ 
if $H \cap gHg^{-1}$ is finite for every $g \in G \setminus H$. 
Hence the condition (2) above can be replaced by the condition that every element of $\H$ is almost malnormal in $G$. 

\begin{lmm}\label{almal}
Let $G$ have a convergence action on a compact metrizable space $X$. 
Then the collection of maximal parabolic subgroups is almost malnormal in $G$. 
\end{lmm}

\begin{proof}
Let $H$ be a maximal parabolic subgroup fixing a parabolic point $p$. 
We have $\Fix(H)=\{p\}$ and $\Stab_G(p)=H$. 
Suppose that $g \in G$ normalizes $H$. 
Then we have $\Fix(H)=\{g(p)\}$ and hence $g$ belongs to $H$. 

Take two different maximal parabolic subgroups $H_1$ and $H_2$ fixing parabolic points $p_1$ and $p_2$, respectively. 
Then the intersection $H_1 \cap H_2$ contains no loxodromic elements and 
fixes two different points $p_1$ and $p_2$. 
Hence $H_1 \cap H_2$ is finite by \cite[Lemma 3B]{Tuk98}. 
\end{proof}

\begin{lmm}\label{order'}
The relation $\to$ defines a partial order 
on the set of almost malnormal and conjugacy invariant collections 
of infinite subgroups of $G$. 
\end{lmm}

\begin{proof}
The relation $\to$ is obviously reflexive and transitive. 
We show that it is antisymmetric. 
Let $\H$ and $\K$ be almost malnormal and 
conjugacy invariant collections of infinite subgroups of $G$ such that $\H \to \K$ 
and $\K \to \H$. We prove that $\H=\K$. 
It suffices to show that $\H \subset \K$. 
Let $H$ be an element of $\H$. 
Then there exists an element $K$ of $\K$ containing $H$. 
Since $\K \to \H$ by the assumption, there exists an element $H'$ of $\H$ containing $K$. 
We have $H \subset K \subset H'$. 
Since $\H$ is almost malnormal in $G$ and $H$ is infinite, this implies that $H=K=H'$.  Therefore $H$ belongs to $\K$. 
\end{proof}

Proposition \ref{order} immediately follows from Lemma \ref{almal} and Lemma \ref{order'}. 

When we consider sets of representatives of conjugacy classes 
of almost malnormal and conjugacy invariant collections
of infinite subgroups of $G$, 
the following is convenient: 
\begin{lmm}\label{rep}
Let $\K$ and $\H$ be almost malnormal and conjugacy invariant collections 
of infinite subgroups of $G$ such that $\K \to \H$. 
If $\{ H_\lambda ~|~ \lambda \in \Lambda\}$ is a set of representatives of conjugacy classes of $\H$, then the following hold: 
\begin{enumerate}
\setlength{\itemsep}{0mm}
\item[(1)] if we take a set $\{ K_{\lambda, \mu} ~|~ \mu \in M_{\lambda}\}$ of representatives of $H_\lambda$-conjugacy classes of $\K_{H_\lambda}$ for each $\lambda \in \Lambda$, then the set $\{ K_{\lambda, \mu} ~|~ \lambda \in \Lambda, \mu \in M_{\lambda}\}$ is a set of representatives of conjugacy classes of $\K$. 
\item[(2)] if the set of conjugacy classes of $\K$ is finite, then the set of $H_\lambda$-conjugacy classes of $\K_{H_\lambda}$ is finite for each $\lambda \in \Lambda$. 
\end{enumerate}
\end{lmm}

\begin{proof} 
(1) Suppose that there exist $\lambda, \lambda' \in \Lambda$, $\mu \in M_\lambda$ and $\mu' \in M_{\lambda'}$ such that $gK_{\lambda,\mu}g^{-1}=K_{\lambda',\mu'}$ for some $g\in G$. 
Then the intersection $gH_\lambda g^{-1} \cap H_{\lambda'}$ is infinite. 
Since $\H$ is almost malnormal in $G$, this implies that $\lambda=\lambda'$ and $g \in H_\lambda$. 
Since $\K_{H_\lambda}$ is also almost malnormal in $H_\lambda$ for every $\lambda \in \Lambda$, we have $\mu=\mu'$ and $g \in K_{\lambda,\mu}$.

On the other hand, since $\K \to \H$ and $\{ H_\lambda ~|~ \lambda \in \Lambda\}$ is a set of representatives of conjugacy classes of $\H$, for every $K\in \K$, there exist $g\in G$ and $\lambda \in \Lambda$ such that $gKg^{-1} \subset H_\lambda$. 
Hence there exist $h \in H_\lambda$ and $\mu \in M_\lambda$ such that $hgK(hg)^{-1}=K_{\lambda, \mu}$. 

\noindent
(2) It follows from (1) that if the set of $H_\lambda$-conjugacy classes of $\K_{H_\lambda}$ is infinite for some $\lambda \in \Lambda$, then the set of conjugacy classes of $\K$ is infinite. 
\end{proof}

\section{The universal relatively hyperbolic structure on a group}

B. Bowditch \cite[Theorem 7.11]{Bow12} characterized relatively hyperbolic structures on a hyperbolic group among conjugacy invariant collections of infinite subgroups. 
If $G$ is a hyperbolic group, then every relatively hyperbolic structure is a blow-down of the empty collection $\emptyset$. 
Taking this into account, we introduce the following notion. 

\begin{dfn}\label{def-universal}
A relatively hyperbolic structure $\K$ on $G$ is called a universal relatively hyperbolic structure on $G$ if every relatively hyperbolic structure on $G$ is a blow-down of $\K$. 
\end{dfn}

Since a universal relatively hyperbolic structure on $G$ is 
the greatest structure with respect to the order $\to$, it is unique if it exists. 
The universal relatively hyperbolic structure is characterized in Corollary \ref{char-universal} for the case of finitely generated groups and we have several finitely generated groups with the universal relatively hyperbolic structure (see Remark (II) in Section \ref{A common}). 

Now we state a characterization of relatively hyperbolic structures on a countable group with the universal relatively hyperbolic structure as follows. This is a consequence of \cite[Theorem 1.1]{Yan11} (refer to Lemma \ref{rep}). 
\begin{prp}\label{rhs}
Let $G$ have the universal relatively hyperbolic structure $\K$ and 
let $\H$ be a conjugacy invariant collection of infinite subgroups of $G$. 
Then $\H$ is a relatively hyperbolic structure on $G$ if and only if the following are satisfied:
\begin{enumerate}
\setlength{\itemsep}{0mm}
\item[(0)] $\K \to \H$;
\item[(1)] $\H$ is almost malnormal in $G$; 
\item[(2)] $\H$ has only finitely many conjugacy classes; 
\item[(3)] every element of $\H$ is quasiconvex relative to $\K$ in $G$.
\end{enumerate}
\end{prp}

\begin{rem}\label{rem-rhs}
We make remarks about the conditions (0) and (3) in Proposition \ref{rhs}. 
\begin{enumerate}
\setlength{\itemsep}{0mm}
\item[(I)] Proposition \ref{rhs} extends Bowditch's characterization of relatively hyperbolic structures on a hyperbolic group \cite[Theorem 7.11]{Bow12}. Indeed, when $G$ is a hyperbolic group, the universal relatively hyperbolic structure of $G$ is the empty collection $\emptyset$ and thus the condition (0) in Proposition \ref{rhs} can be omitted. 
However, we cannot omit the condition (0) in Proposition \ref{rhs} in general
(see Proposition \ref{pA}). 
\item[(II)] If $G$ is a finitely generated group, then the condition (3) in Proposition \ref{rhs} can be replaced the following condition: 
\begin{enumerate}
\setlength{\itemsep}{0mm}
\item[(3)'] every element of $\H$ is finitely generated and undistorted in $G$. 
\end{enumerate}

Indeed if $\H$ is a relatively hyperbolic structure on $G$, then every element of $\H$ is finitely generated by \cite[Proposition 2.29]{Osi06a} and it is undistorted in $G$ by \cite[Lemma 5.4]{Osi06a}. 
On the other hand, if we suppose that the condition (3)' holds, then it follows from \cite[Theorem 1.5]{Hru10} that every element of $\H$ is quasiconvex relative to $\K$ in $G$. 
\item[(III)] There exist finitely generated groups $G$ with the universal relatively hyperbolic structure $\K$ such that every finitely generated subgroup of $G$ is quasiconvex relative to $\K$ in $G$ (see Remark (III) in Section \ref{A common}). 
For such groups, we can replace the condition (3) in Proposition \ref{rhs} by the following condition:
\begin{enumerate}
\setlength{\itemsep}{0mm}
\item[(3)''] every element of $\H$ is finitely generated. 
\end{enumerate}
We do not know whether (3) in Proposition \ref{rhs} can be replaced by (3)'' whenever $G$ is a finitely generated group with the universal relatively hyperbolic structure.
\end{enumerate}
\end{rem}

\section{Relative quasiconvexity under blowing up and down a relatively hyperbolic structure}
We characterize how relative quasiconvexity for subgroups varies when we blow up and down a relatively hyperbolic structure as follows:
\begin{prp}\label{qc-updown}
Let $\K$ and $\H$ be two relatively hyperbolic structures on $G$ such that $\K\to \H$. 
Then for every subgroup $L$ of $G$, the following conditions (i), (ii), (iii) are equivalent: 
\begin{enumerate}
\setlength{\itemsep}{0mm} 
\item[(i)] $L$ is quasiconvex relative to $\K$ in $G$; 
\item[(ii)] $L$ is quasiconvex relative to $\H$ in $G$ 
and $L\cap H$ is quasiconvex relative to $\K$ in $G$ 
for every $H \in \H$; 
\item[(iii)] $L$ is quasiconvex relative to $\H$ in $G$ 
and $L\cap H$ is quasiconvex relative to $\K_H$ in $H$ 
for every $H \in \H$.
\end{enumerate}
\end{prp}

The equivalence of the conditions (i) and (ii) follows from \cite[Theorem 1.3]{Yan11} 
and Lemma \ref{rep}. 
The equivalence of the conditions (ii) and (iii) is implied by the following: 
\begin{lmm}\label{relqc-rest}
Let $\K$ and $\H$ be relatively hyperbolic structures on $G$ such that $\K \to \H$. 
Fix an element $H$ of $\H$. 
Then for every subgroup $L$ of $H$, the following conditions (i), (ii) are equivalent:
\begin{enumerate}
\setlength{\itemsep}{0mm}
\item[(i)] $L$ is quasiconvex relative to $\K$ in $G$;
\item[(ii)] $L$ is quasiconvex relative to $\K_H$ in $H$.
\end{enumerate}
\end{lmm}

\begin{proof}
First we remark that $\K_H$ is a relatively hyperbolic structure on $H$ by 
\cite[Corollary 3.4]{Yan11}. 
It follows from the definition that $\K_H \subset H\wedge \K$. 
Let $K$ be an element of $\K$ such that $H \cap K$ is infinite. 
Since $\K \to \H$, there exists $H' \in \H$ containing $K$. 
Then $H \cap H'$ is infinite. 
Since $\H$ is a relatively hyperbolic structure on $G$, this implies that $H=H'$. 
Hence $K$ is contained in $H$ and it belongs to $\K_H$. 
Thus we showed that $\K_H=H\wedge \K$. 
Hence the assertion follows from \cite[Corollary 9.3]{Hru10}. 
\end{proof}

In order to give a criterion for 
two relatively hyperbolic structures on a countable group 
to give the same set of relatively quasiconvex subgroups, 
we show the following: 
\begin{prp}\label{nqc}
Let $G$ be endowed with a relatively hyperbolic structure $\K$. 
Then every subgroup of $G$ is quasiconvex relative to $\K$ in $G$ if and only if either $G$ is virtually cyclic or $\K$ is trivial. 
\end{prp}
\begin{proof}
The `if' part follows from Remark (I) in Section \ref{Relatively}. 
We prove the `only if' part. 
Suppose that $G$ is not virtually cyclic and that $\K$ is proper. 
Then it follows from \cite[Corollary 3.2]{M-O-Y4} 
that there exists a subgroup $L$ of $G$ which is a free group of rank two 
and strongly quasiconvex relative to $\K$ in $G$. 
Then $L$ has a subgroup $L'$ which is not finitely generated. 
Assume that $L'$ is quasiconvex relative to $\K$ in $G$. 
Then $L' \wedge \K$ is a relatively hyperbolic structure on $L'$ 
by \cite[Theorem 1.2 (1)]{Hru10}. 
Since $L$ is strongly quasiconvex relative to $\K$ in $G$, 
we have $L \wedge \K=\emptyset$ and hence $L' \wedge \K=\emptyset$. 
This implies that $L'$ is a hyperbolic group, which contradicts the fact 
that hyperbolic groups are finitely generated. 
Therefore $L'$ is not quasiconvex relative to $\K$ in $G$. 
\end{proof}

For a relatively hyperbolic structure $\H$ on $G$, 
note that $\H_{\ne}$ is a relatively hyperbolic structure on $G$ 
(see \cite[Theorem 2.40]{Osi06a}).
If two relatively hyperbolic structures $\K$ and $\H$ on $G$
satisfy $\K_{\ne}=\H_{\ne}$, then we have $\RQC(G,\K)=\RQC(G,\H)$
(\cite[Corollary 1.3]{MP08}). On the other hand we have the following:
\begin{crl}\label{qc-updown1}
Let $\K$ and $\H$ be two relatively hyperbolic structures on $G$. 
If $\K \to \H$ and $\RQC(G,\K)=\RQC(G,\H)$, then we have $\K_{\ne}=\H_{\ne}$. 
\end{crl}
\noindent
For the case of finitely generated groups, 
this is refined as in Theorem \ref{qc-updown2}. 

\begin{proof}[Proof of Corollary \ref{qc-updown1}]
We suppose that $\K \to \H$ and $\RQC(G,\K)=\RQC(G,\H)$. 
It follows from the equivalence between the conditions (i) and (iii) in Proposition \ref{qc-updown} 
that for every $H \in \H_{\ne}$, every subgroup of $H$ is quasiconvex relative to $\K_H$ in $H$. 
Since $H$ is not virtually cyclic, we have $\K_{H}=\{H\}$ by Proposition \ref{nqc}. 
This implies that $\H_{\ne} \to \K_{\ne}$ and hence $\K_{\ne}=\H_{\ne}$. 
\end{proof}

\section{Cardinality of the set of relatively hyperbolic structures} 
First we consider virtually infinite cyclic groups. 

\begin{prp}\label{virtinfcyclic}
The group $G$ is virtually infinite cyclic if and only if $\RHS(G)=\{\emptyset, \{ G \}\}$. 
\end{prp}

In order to prove this, 
we need the following:
\begin{lmm}\label{inf-ind}
Let $\H$ be an almost malnormal and conjugacy invariant collection of infinite subgroups of $G$. 
If $\H$ is not equal to $\{G\}$, then $\H$ consists of infinite index subgroups of $G$. 
\end{lmm}
\begin{proof}
Assume that an element $H\in \H$ is a finite index subgroup of $G$. 
Then $H$ contains a finite index normal subgroup $H'$ of $G$. 
Since $G$ is infinite, $H'$ is also infinite. 
For every $g \in G$, we have $H \cap gHg^{-1} \supset H'$. 
In particular for every $g \in G \setminus H$, 
$H \cap gHg^{-1}$ is infinite. 
This contradicts the assumption that $\H$ is almost malnormal in $G$. 
\end{proof}

\begin{proof}[Proof of Proposition \ref{virtinfcyclic}]
Since $G$ is infinite and hyperbolic if and only if 
$\{ G \}$ and $\emptyset$ are elements of $\RHS(G)$, 
we suppose that $G$ is infinite and hyperbolic. 

If $G$ is virtually infinite cyclic with a relatively hyperbolic structure $\H \neq \emptyset$, then every element of $\H$ is a finite index subgroup of $G$. 
Hence $\H$ is trivial by Lemma \ref{inf-ind}. 
On the other hand, if $G$ is not virtually infinite cyclic then it follows from \cite[Corollary 1.7]{Osi06b} that $G$ has a virtually infinite cyclic subgroup $H$ which is hyperbolically embedded into $G$ relative to the empty collection $\emptyset$. 
\end{proof}

Second we consider countable groups which are not virtually cyclic. 
\begin{thm}\label{inf-rhs}
Suppose that $G$ is not virtually cyclic and admits a proper relatively hyperbolic structure $\H$. 
Then there exists a sequence $(\H_n)_{n \in \N \cup \{0\}}$ of proper relatively hyperbolic structures on $G$ satisfying the following: 
\begin{enumerate}
\setlength{\itemsep}{0mm}
\item[(1)] $\H_0=\H$ and $\H_{n} \subsetneqq \H_{n+1}$ for every $n \in \N \cup \{0\}$; 
\item[(2)] if $i < j$, then $\RQC(G, \H_i)$ is a proper subset of $\RQC(G, \H_j)$. 
\end{enumerate}
On the other hand there exists a sequence $(\K_n)_{n \in \N \cup \{0\}}$ of proper relatively hyperbolic structures on $G$ satisfying the following: 
\begin{enumerate}
\setlength{\itemsep}{0mm}
\item[(1)] $\K_0=\H$ and $\K_{n} \subsetneqq \K_{n+1}$ for every $n \in \N \cup \{0\}$; 
\item[(2)] $\RQC(G, \K_n)=\RQC(G, \H)$ for every $n\in \N \cup \{0\}$. 
\end{enumerate}
In particular $\RHS(G)/\Out(G)$ has infinitely many elements. 
\end{thm}

\begin{proof}
We construct $\H_n$ inductively. 
We put $\H_0=\H$. 
Suppose that $\H_n$ is constructed. 
Since $\H_n$ is proper, it follows from \cite[Theorem 1.2]{M-O-Y4} that there exists a finitely generated and virtually non-abelian free subgroup $V$ of $G$ which is hyperbolically embedded into $G$ relative to $\H_n$. 
We put $\H_{n+1}=\H_n \cup \{ P ~|~ P=gVg^{-1}$ for some $g \in G\}$. 
If $i < j$, then it follows from \cite[Theorem 1.1(1) and Corollary 1.3]{MP08} 
that $\RQC(G, \H_i)$ is a proper subset of $\RQC(G, \H_j)$. 

Also we have a sequence $(\K_n)_{n\in \N \cup \{0\}}$ with the desired properties by using \cite[Corollary 1.7]{Osi06b} instead of \cite[Theorem 1.2]{M-O-Y4} on the above argument. 
\end{proof}

For finitely generated groups, we have the following. 
\begin{prp}\label{fingen}
If $G$ is a finitely generated group, 
then $\RHS(G)$ is countable. 
\end{prp}
\begin{proof}
For a set $A$, we denote by $F(A)$ the set of all finite subsets of $A$. 
We suppose that $G$ is a finitely generated group and 
construct an injective map $i:\RHS(G) \to F(F(G))$ as follows. 

Let $\H$ be a relatively hyperbolic structure on $G$. 
We take a finite set $\{ H_i ~|~ i \in \{1, \ldots, n\} \}$ of representatives  of conjugacy classes of $\H$. 
Since $G$ is finitely generated, 
$\H$ consists of finitely generated subgroups of $G$ by 
\cite[Proposition 2.29]{Osi06a}. 
Hence we can choose a finite generating set $S_i$ of $H_i$ 
for each $i\in\{1, \ldots, n\}$. 
We define $i(\H)=\{ S_i ~|~ i \in \{1, \ldots, n\} \}$. 

Since $G$ is finitely generated, it is countable and hence $F(F(G))$ is also countable. 
Therefore $\RHS(G)$ is also countable. 
\end{proof}

On the other hand we have the following: 
\begin{prp}
Let $G_k$ be an infinite countable group for each positive integer $k$ and let $G=*_{k\in \N}G_k$. 
Then $\RHS(G)$ is uncountable. 
Moreover if we suppose that for every $k \in \N$, $G_k$ is freely indecomposable 
and not isomorphic to $G_j$ for every $j \in \N \setminus \{k\}$, 
then $\RHS(G) \slash \Out(G)$ is also uncountable. 
\end{prp}

\begin{proof}
We denote by $\{0,1\}^\N$ the set of all maps from $\N$ to $\{0,1\}$. 
For each $\sigma\in\{0,1\}^\N$, we consider a relatively hyperbolic structure $\H_\sigma$ represented by $\{*_{\sigma(k)=0}G_k, *_{\sigma(k)=1}G_k\}$. 
Since we have $\H_{\sigma_1}\neq\H_{\sigma_2}$ for $\sigma_1, \sigma_2 \in \{0,1\}^\N$ such that $\sigma_1\neq \sigma_2$, the former assertion follows from the fact that $\{0,1\}^\N$ is uncountable. 
Moreover if we suppose that for every $k \in \N$, $G_k$ is freely indecomposable 
and not isomorphic to $G_j$ for every $j \in \N \setminus \{k\}$, then it follows from Kurosh subgroup theorem (see for example \cite[Chapter IV, Theorem 1.10]{L-S77}) that no automorphisms of $G$ transform $\H_{\sigma_1}$ to $\H_{\sigma_2}$. 
\end{proof}

\section{A common blow-up of two relatively hyperbolic structures}\label{A common}
In the case of finitely generated groups, we have the following: 

\begin{prp}\label{fiberprod}
Let $G$ be a finitely generated group. 
If $\H$ and $\K$ are relatively hyperbolic structures on $G$, 
then $\H\wedge \K$ is also a relatively hyperbolic structure on $G$. 
\end{prp}

In the above, $\H\wedge \K$ is a common blow-up of two relatively hyperbolic structures $\H$ and $\K$. 
In particular, the partially ordered set $(\RHS(G), \to)$ is a directed set. 

In order to prove Proposition \ref{fiberprod}, 
we need the following:
\begin{lmm}\label{rest}
Let $G$ be a finitely generated group with two relatively hyperbolic structures $\H$ and $\K$.
Then for every $H\in \H$, $H\wedge \K$ is a relatively hyperbolic structure on $H$.
\end{lmm}

\begin{proof}
Since $H$ is finitely generated by \cite[Proposition 2.29]{Osi06a}, it follows from \cite[Lemma 5.4]{Osi06a} that $H$ is undistorted in $G$. 
Therefore $H\wedge \K$ is a relatively hyperbolic structure on $H$ by \cite[Theorem 1.8]{D-S05b} (see also \cite[Theorem 1.5 and Theorem 9.1]{Hru10}).
\end{proof}

\begin{proof}[Proof of Proposition \ref{fiberprod}]
We denote the conjugacy classes of $\H$ by $\H_1, \ldots, \H_n$. 
Then we have $\H=\bigsqcup_{i=1}^n \H_i$ and $\bigsqcup_{i=1}^{n} (\H_i \wedge \K)=\H \wedge \K$. 
We show that for each $l \in \{1, \ldots, n+1\}$, $(\bigsqcup_{i=1}^{l-1} (\H_i \wedge \K)) \cup (\bigsqcup_{i=l}^{n} \H_i)$ is a relatively hyperbolic structure on $G$. 
The proof is done by induction on $l$. 
When $l=1$, the assertion obviously holds. 
Suppose that $m \in  \{1, \ldots, n\}$ and that $(\bigsqcup_{i=1}^{m-1} (\H_i \wedge \K)) \cup (\bigsqcup_{i=m}^{n} \H_i)$ is a relatively hyperbolic structure on $G$. 
Each element $H$ of $\H_m$ is hyperbolic relative to $H \wedge \K$ by Lemma \ref{rest}. 
Since we have $\bigsqcup_{H \in \H_m} (H \wedge \K)=\H_m \wedge \K$, 
it follows from \cite[Corollary1.14]{D-S05b} that $G$ is hyperbolic relative to 
$(\bigsqcup_{i=1}^{m} (\H_i \wedge \K)) \cup (\bigsqcup_{i=m+1}^{n} \H_i)$.  
\end{proof}

Now we can characterize the universal relatively hyperbolic structure for the case of  finitely generated groups. 
\begin{crl}\label{char-universal}
Let $G$ be a finitely generated group and 
let $\K$ be a relatively hyperbolic structure on $G$. 
Then the following are equivalent: 
\begin{enumerate}
\item[(i)] $\K$ is the universal relatively hyperbolic structure on $G$;
\item[(ii)] every element of $\K$ has no proper relatively hyperbolic structures;
\item[(iii)] no relatively hyperbolic structure on $G$ other than $\K$ is a blow-up of $\K$.
\end{enumerate}
\end{crl}

\begin{proof}
The implication $(i)\Rightarrow  (iii)$ follows from Proposition \ref{order}. 

We prove that $(iii)\Rightarrow  (i)$.
Suppose that no relatively hyperbolic structure on $G$ other than $\K$ is a blow-up of $\K$. 
Let $\H$ be an arbitrary relatively hyperbolic structure on $G$. 
Then it follows from Proposition \ref{fiberprod} that 
$\H\wedge \K$ is a relatively hyperbolic structure on $G$ 
such that $\H\wedge \K \to \H$ and $\H\wedge \K \to \K$. 
We have $\H \wedge \K=\K$ by the assumption and hence $\H$ is a blow-down of $\K$. 

Next we prove that $(ii)\Rightarrow  (iii)$.
Suppose that every element of $\K$ has no proper relatively hyperbolic structures. 
Let $\H$ be an arbitrary relatively hyperbolic structure on $G$ with $\H \to \K$. 
We prove that $\H=\K$. 
Let $K$ be an arbitrary element of $\K$. 
It follows from \cite[Corollary 3.4]{Yan11} that $\H_K$ 
is a relatively hyperbolic structure on $K$. 
The assumption on $\K$ implies that $\H_K=\{K\}$. 
Hence $K$ is an element of $\H$ and this implies that $\K \to \H$. 
It follows from the assumption on $\K$ that $\H \to \K$ and hence we have $\H=\K$ 
by Proposition \ref{order}. 

Finally we prove that $(iii)\Rightarrow  (ii)$.
Suppose that there exists an element $K$ of $\K$ 
which has a proper relatively hyperbolic structure $\K'$. 
We set $\H=(\K \setminus \{ L ~|~ L=gKg^{-1}$ for some $g \in G\}) 
\cup \{ P ~|~ P=gQg^{-1}$ for some $Q \in \K'$ and $g \in G\}$. 
It follows from \cite[Corollary 1.14]{D-S05b} that $\H$ 
is a relatively hyperbolic structure on $G$. 
It follows from the construction of $\H$ that $\H \neq \K$ and $\H$ is a blow-up of $\K$. 
\end{proof}

\begin{rem}\label{rem-universal}
We make remarks on the universal relatively hyperbolic structure. 
\begin{enumerate}
\setlength{\itemsep}{0mm}

\item[(I)] Since the action of the outer automorphism group $\Out(G)$ on $\RHS(G)$ preserves the order $\to$, the universal relatively hyperbolic structure on $G$ is invariant under the action of $\Out(G)$ (see also \cite[Lemma 4.23(4)]{D-S08}). 
For finitely generated groups, 
the property of having no proper relatively hyperbolic structures 
is a quasi-isometric invariant by \cite[Theorem 1.2]{Dru09}. 
In view of \cite[Theorem 4.8]{B-D-M09}, Corollary \ref{char-universal} implies that for finitely generated groups the existence of the universal relatively hyperbolic structure is invariant under quasi-isometry. 
For a finitely generated group $G$ with the universal relatively hyperbolic structure, relationship between splitting of $G$ and $\Out(G)$ is described in \cite[Theorem 1.12]{D-S08}. 

\item[(II)] If $G$ is infinite and has no proper relatively hyperbolic structures, 
then the trivial relatively hyperbolic structure $\K=\{G\}$ 
is a unique relatively hyperbolic structure on $G$ and hence it is universal. 
The examples of such groups are $\Z^n\ (n \ge 2)$, $\SL(n,\Z)\ (n \ge 3)$ 
and the mapping class group of an orientable surface of genus $g$ with $p$ punctures, 
where $3g+p \ge 5$ (see \cite[Theorem 1.2 and p.557]{B-D-M09} 
and \cite[Section 8]{K-N04} for details and other examples).
There exists a criterion for countable groups to have no proper relatively hyperbolic structures 
(see \cite[Theorem 1]{K-N04} together with \cite[Definition 1]{Bow12}, 
and also \cite[Theorem 2]{A-A-S07}). 

\indent
On the other hand, Corollary \ref{char-universal} enables us to recognize that each of the following finitely generated groups has the universal relatively hyperbolic structure $\K$, which is proper: 
\begin{itemize}
\setlength{\itemsep}{0mm}
\item each hyperbolic group with $\K=\emptyset$; 
\item each geometrically finite Kleinian group with the collection $\K$ of all maximal parabolic subgroups that are not virtually infinite cyclic (note that every maximal parabolic subgroup of a Kleinian group is virtually abelian (see for example \cite[Proposition 2.2]{M-T98})); 
\item each free product $A*B$ with the collection $\K$ of all conjugates of $A$ and $B$, where $A$ and $B$ are finitely generated groups having no proper relatively hyperbolic structures;
\item each one-relator product $A*B \slash \langle\langle r^m \rangle\rangle$ with the collection $\K$ of all conjugates of $A$ and $B$, where $A$ and $B$ are finitely generated groups having no proper relatively hyperbolic structures, $r$ is a cyclically reduced word of length at least 2 and $m \ge 6$ (see \cite[Theorem 4.1]{MP-W11b}); 
\item each limit group with the collection $\K$ of all maximal abelian non-cyclic subgroups (\cite[Theorem 0.3]{Dah03}). 
\end{itemize}

\item [(III)] The following are examples of finitely generated groups $G$ 
with the universal relatively hyperbolic structure $\K$ such that 
every finitely generated subgroup of $G$ is quasiconvex relative to $\K$ in $G$ 
(refer to the above (II)): 
\begin{itemize}
\setlength{\itemsep}{0mm}
\item each finitely generated free group (\cite[Section 2]{Sho91}) and the fundamental group of each closed hyperbolic surface (\cite[Proposition 2]{Pit93});
\item geometrically finite Kleinian groups of the second kind acting on $\mathbb{H}^3$ (\cite[Proposition 7.1]{Mor84}); 
\item each free product $A*B$ where $A$ and $B$ are finitely generated groups having no proper relatively hyperbolic structures;
\item each one-relator product $A*B \slash \langle\langle r^m \rangle\rangle$ where $A$ and $B$ are finitely generated groups having no proper relatively hyperbolic structures, $r$ is a cyclically reduced word of length $|r| \ge 2$ and $m > 3|r|$ (\cite[Theorem 1.7]{MP-W11b}); 
\item each limit group (\cite[Proposition 4.6]{Dah03}). 
\end{itemize}

\item [(IV)] There exist finitely generated groups which do not have the universal relatively hyperbolic structure. 
An example of such groups is the so-called Dunwoody's inaccessible group (see \cite[Section 6]{B-D-M09}). 
Note that Dunwoody's inaccessible group is finitely generated, not finitely presentable and has torsions. 
We can also obtain torsion-free countable non-finitely generated groups without the universal relatively hyperbolic structure (see Proposition \ref{freeprod}).  
However, it is unknown whether every finitely presented (resp. torsion-free finitely generated) group has the universal relatively hyperbolic structure (see \cite[Question 1.5]{B-D-M09}). 
\end{enumerate}
\end{rem}

\section{Relatively hyperbolic structures with the same set of relatively quasiconvex subgroups}

We determine when two relatively hyperbolic structures 
have the same collection of relatively quasiconvex subgroups. 

\begin{thm}\label{qc-updown2}
Let $G$ be a finitely generated group and let 
$\K$ and $\H$ be two relatively hyperbolic structures on $G$. 
Then we have the following:
\begin{enumerate}
\setlength{\itemsep}{0mm} 
\item[(1)] $\RQC(G,\K)\cap\RQC(G,\H)=\RQC(G,\K\wedge\H)$.
\item[(2)] The following conditions are equivalent:
\begin{enumerate}
\setlength{\itemsep}{0mm}
\item[(i)] $\RQC(G,\K)=\RQC(G,\H)$; 
\item[(ii)] $\RQC(G,\K)=\RQC(G,\H)=\RQC(G,\K\wedge\H)$;
\item[(iii)] $\K_{\ne}= \H_{\ne}$. 
\end{enumerate}
\end{enumerate}
\end{thm}
\noindent
We remark that 
the equivalence between (i) and (iii) in (2) follows from \cite[Corollary 1.3]{MP08} 
when $\K$ is a subcollection of $\H$. 

\begin{proof}
(1) Let $L$ be a subgroup of $G$. 
First we suppose that $L$ is quasiconvex relative to $\K\wedge\H$ in $G$. 
Since both $\K$ and $\H$ are blow-downs of $\K\wedge\H$, it follows from Proposition \ref{qc-updown} (i) $\Rightarrow$ (ii) that $L$ is quasiconvex relative to $\K$ in $G$ and quasiconvex relative to $\H$ in $G$. 

Next we suppose that $L$ is quasiconvex relative to $\K$ in $G$ and quasiconvex relative to $\H$ in $G$. 
Let $H$ be an element of $\H$. 
Since $G$ is finitely generated, $H$ is undistorted in $G$ by \cite[Lemma 5.4]{Osi06a}. 
It follows from \cite[Theorem 1.5]{Hru10} that $H$ is quasiconvex relative to $\K$ in $G$. 
Therefore we have $L \cap H$ is also quasiconvex relative to $\K$ in $G$ by \cite[Theorem 1.2 (2)]{Hru10}. 
Since $H$ and $L\cap H$ are quasiconvex relative to $\K$ in $G$ and $L\cap H$ is a subgroup of $H$, 
$L \cap H$ is quasiconvex relative to $H \wedge \K$ in $H$ by \cite[Corollary 9.3]{Hru10}
(see also Lemma \ref{relqc-rest}). 
Since we have $H \wedge \K=(\K \wedge \H)_H$, $L \cap H$ is quasiconvex relative to $(\K \wedge \H)_H$ in $H$. 
Hence it follows from Proposition \ref{qc-updown} (iii) $\Rightarrow$ (i) that $L$ is quasiconvex relative to $\K\wedge\H$ in $G$. 

(2) Since we have $\RQC(G,\K) \cap \RQC(G,\H)=\RQC(G, \K\wedge \H)$, the 
conditions (i) and (ii) are equivalent. 

The condition (iii) implies the condition (i) 
by \cite[Corollary 1.3]{MP08}.

Finally we prove that the condition (ii) implies the condition (iii). 
Since we have $\K \wedge \H \to \H$ and $\K \wedge \H \to \K$, 
it follows from Corollary \ref{qc-updown1} that $\H_{\ne}=(\K \wedge \H)_{\ne}=\K_{\ne}$.  
\end{proof}

\appendix
\def\thesection{Appendix \Alph{section}}
\section{Torsion-free countable groups without the universal relatively hyperbolic structure}
\def\thesection{\Alph{section}}

We give examples of torsion-free countable groups without the universal relatively hyperbolic structure. 
Note that these are not finitely generated. 

\begin{prp}\label{freeprod}
Let $G_l$ be an infinite countable group for each $l \in \N$ 
and put $G=*_{l\in \N}G_l$.
Then $G$ has no universal relatively hyperbolic structures. 
In particular, infinite countably generated free groups have 
no universal relatively hyperbolic structures. 
\end{prp}

\begin{proof}
Assume that $G$ has the universal relatively hyperbolic structure $\K$. 
For each $m \in \N$, we put $A_m=*_{l=1}^mG_l$ and $Z_m=*_{l=m+1}^{\infty}G_l$.  We have $G=A_m * Z_m$ and hence $G$ has a relatively hyperbolic structure $\H_{m}$ represented by $\{A_m, Z_m\}$. 
Since $\K$ is universal, it is a blow-up of $\H_m$ for every $m \in \N$. 
Hence it follows from \cite[Corollary 3.4]{Yan11} 
that $\K_{Z_m}$ is a relatively hyperbolic structure on $Z_m$ for each $m \in \N$. 
Since $Z_m$ is not finitely generated, it is not hyperbolic. 
Therefore $\K_{Z_m} \neq \emptyset$ for every $m \in \N$. 
Since $\K$ has only finitely many conjugacy classes, there exists a conjugacy class $\K'$ of $\K$ such that $\K' \cap \K_{Z_m} \neq \emptyset$ for every $m\in \N$. 
Let $K'$ be a representative of $\K'$ and let $k'$ be a nontrivial element of $K'$. 
There exists $n \in \N$ such that $k'$ belongs to $A_n$. 
Since $\K' \cap \K_{Z_n} \neq \emptyset$, there exist an element $g$ of $G$ such that $gk'g^{-1}$ belongs to $Z_n$. 
Thus we have $gA_ng^{-1} \cap Z_n \neq \{1\}$. 
However, this contradicts the fact that we have $G=A_n * Z_n$. 
\end{proof}

\def\thesection{Appendix \Alph{section}}
\section{A remark on Proposition \ref{rhs}}
\def\thesection{\Alph{section}}

We give a finitely generated group $G$ with a conjugacy invariant collection of infinite subgroups $\H$ such that $\H$ satisfies the conditions (1), (2) in Proposition \ref{rhs} and the condition (3)' in Remark \ref{rem-rhs} (II) and it is not a relatively hyperbolic structure on $G$. 

Let $G$ be the mapping class group of an orientable surface of genus $g$ with $p$ punctures, where $3g+p \ge 5$. 
Then $G$ has no proper relatively hyperbolic structures and the trivial relatively hyperbolic structure $\K=\{G\}$ is the universal relatively hyperbolic structure on $G$ (see Remark \ref{rem-universal} (II)). 
Let $P$ be a subgroup of $G$ generated by a pseudo-Anosov element and let $V_G(P)$ be the virtual normalizer of $P$ in $G$, that is, $V_G(P)=\{ g \in G ~|~ [P \colon P \cap gPg^{-1}] < \infty$ and $[gPg^{-1} \colon P \cap gPg^{-1}] < \infty\}$. 
We denote by $\H$ the set of all conjugates of $V_G(P)$ in $G$. 
Since every pseudo-Anosov element of $G$ is of infinite order, $\H$ is a conjugacy invariant collection of infinite subgroups of $G$. 

\begin{prp}\label{pA}
$\H$ satisfies the conditions (1), (2) in Proposition \ref{rhs} and the condition (3)' in Remark \ref{rem-rhs} (II). 
\end{prp}

In order to prove this, 
we need the following:
\begin{lmm}\label{vn}
Let $L$ be a group and let $M$ be an infinite subgroup of $L$. 
Suppose that the virtual normalizer $V_L(M)$ of $M$ is virtually infinite cyclic. Then $V_L(M)$ is an almost malnormal subgroup of $L$. 
\end{lmm}

\begin{proof}
It follows from the assumption that $M$ is a finite index subgroup of $V_L(M)$.  It follows that for every $l \in L$, we have $[V_L(M) \cap lV_L(M)l^{-1} \colon M \cap lMl^{-1}] < \infty$. 

Suppose that $l$ belongs to $V_L(V_L(M))$. 
Then we have $[V_L(M) \colon V_L(M) \cap lV_L(M)l^{-1}] < \infty$ and $[lV_L(M)l^{-1} \colon V_L(M) \cap lV_L(M)l^{-1}] < \infty$. 
It follows that we have $[V_L(M) \colon M \cap lMl^{-1}] < \infty$ and $[lV_L(M)l^{-1} \colon M \cap lMl^{-1}] < \infty$. 
This implies that $l$ belongs to $V_L(M)$. 
Thus we have $V_L(V_L(M))=V_L(M)$. 

Now suppose that $l$ belongs to $L \setminus V_L(M)$. 
Then neither $l$ nor $l^{-1}$ belongs to $V_L(M)$. 
Hence the equality $V_L(V_L(M))=V_L(M)$ implies that $V_L(M) \cap lV_L(M)l^{-1}$ is an infinite index subgroup of $V_L(M)$. 
Since $V_L(M)$ is virtually infinite cyclic, $V_L(M) \cap lV_L(M)l^{-1}$ is finite. 
\end{proof}

\begin{proof}[Proof of Proposition \ref{pA}]
It is known that 
$V_G(P)$ is a virtually infinite cyclic subgroup of $G$ 
(see for example \cite[Theorem 3.5]{Mos07})
and hence $P$ is a finite index subgroup of $V_G(P)$. 
Hence the condition (1) holds by Lemma \ref{vn}. 
Since $\H$ consists of all conjugates of a single subgroup $V_G(P)$, the condition (2) also holds. 
Every free abelian subgroup of $G$ is undistorted in $G$ (see \cite[Corollary 5.3 (1)]{Ham06}) and hence $P$ is undistorted in $G$. 
Since $P$ is a finite index subgroup of $V_G(P)$, $V_G(P)$ is also undistorted in $G$. 
Hence the condition (3)' holds. 
\end{proof}

On the other hand, $G$ has no proper relatively hyperbolic structure as we mentioned above and hence $\H$ is not a relatively hyperbolic structure on $G$.


\end{document}